\documentclass[12pt]{article}
\usepackage{geometry} 
\usepackage{array} 
\usepackage{paralist} 
\usepackage{verbatim} 

\usepackage[utf8]{inputenc} 
\usepackage{textcomp} 
\usepackage{graphicx}  
\usepackage{epstopdf} 
\usepackage{flafter}  

\usepackage{amsfonts,latexsym,amsthm,amssymb}
\usepackage{amsmath,amssymb}  
\usepackage{bm}  

\usepackage{fancyhdr} 
\usepackage{sectsty}
\allsectionsfont{\sffamily\mdseries\upshape} 

\numberwithin{equation}{section}

\def\Ref{\ref}
\def\c{\cite}
\def\Ref#1{(\ref{#1})}

\begin{document}
\title{Quantum symmetry algebras of spin systems related to Temperley - Lieb R-matrices}

\author{P. ~P.~Kulish $^1$
\thanks{E-mail address: kulish@pdmi.ras.ru}
\and N. ~Manojlovi\'c $^{2,3}$
\thanks{E-mail address: nmanoj@ualg.pt}
\and Z. ~Nagy $^{2,3}$
\thanks{E-mail address: znagy@ualg.pt} \\
\\
$^1$\textit{St. Petersburg Department of Steklov Mathematical Institute} \\
\textit{Fontanka 27, 191023, St.Petersburg, Russia} \\
\\
$^2$\textit{Departamento de Matem\'atica, F. C. T.,
Universidade do Algarve}\\ 
\textit{Campus de Gambelas, PT-8005-139 Faro, Portugal} \\ 
\\
$^3$\textit{Grupo de F\'{\i}sica Matem\'atica da Universidade de Lisboa} \\
\textit{Av. Prof. Gama Pinto 2, PT-1649-003 Lisboa, Portugal}}

\date{}

\maketitle



\begin{abstract}
A reducible representation of the Temperley-Lieb algebra is constructed on the tensor product of $n$-dimensional spaces. One obtains as a centraliser of this action a quantum algebra 
(a quasi-triangular Hopf algebra) $\mathcal{U}_q$ with a representation ring equivalent to the 
representation ring of the $\mathfrak{sl}_2$ Lie algebra. This algebra $\mathcal{U}_q$ is
the symmetry algebra of the corresponding open spin chain.
\end{abstract} 


\clearpage \newpage

%
%

\section{Introduction}

The development of the quantum inverse scattering method (QISM) 
\c{SklyaninTakhFad,TakhFadI, Fadddeev, KS} as an approach to resolution and construction of quantum integrable systems has lead to the foundations of the theory of quantum groups \c{KR,Sklyanin, Drinfeld,Jimbo}.

The theory of representations of quantum groups is naturally connected to the spectral theory of the integrals of motion of quantum systems. In particular, this connection appeared in the combinatorial approach to the question of completeness of the eigenvectors of the $XXX$ Heisenberg spin chain \c{TakhFadII} 
with the Hamiltonian
\begin{equation}
\label{HXXX}
H_{XXX} = \sum _{n=1}^{N-1} \left( \sigma_n^x \sigma_{n+1}^x + \sigma_n^y \sigma_{n+1}^y  + \sigma_n^z \sigma_{n+1}^z  \right) ,
\end{equation}
where $\sigma_n^{\alpha}, \alpha = x,y,z$ are Pauli matrices.

Three algebras are connected to this system: the Lie algebra 
$\mathfrak{sl}_2$ of rotations, the group algebra 
$\mathbb{C}[S_N]$ of the symmetric group $S_N$ and the infinite 
dimensional algebra $\mathcal{Y}(\mathfrak{sl}_2)$ -- the Yangian, \c{Drinfeld} with the corresponding $R$-matrix $R(\lambda)=\lambda I+\eta \mathcal{P}$, where $\mathcal{P}$ is the $4\times 4$ permutation matrix flipping 
the two factors of $\mathbb{C}^2 \otimes \mathbb{C}^2$.

The Yangian is the dynamical symmetry algebra which contains 
all the dynamical observables of the system. It is important to 
note that the algebras $\mathfrak{sl}_2$ and $\mathbb{C}[S_N]$ are 
related by the Frobenius-Schur-Weyl duality in the representation 
space $\mathcal{H}= \underset {1}{\overset {N}{\otimes}}\mathbb{C}^2$. 
This follows from the fact that $\mathfrak{sl}_2$ and $\mathbb{C}[S_N]$ are each other's centralizers in this representation space. As a consequence, since the Hamiltonian commutes with the global generators of 
$\mathfrak{sl}_2: S^{\alpha} =1/2 \sum _{n=1}^N \sigma_n^{\alpha} , \alpha = x,y,z,$ it is an element of $\mathbb{C}[S_N]$.  
This can also be seen from the expression of $H_{XXX}$ in terms of the permutation 
operators, which are the generators of the symmetric group $S_N$
$$
\sum _{\alpha} \sigma_n^{\alpha} \sigma_{n+1}^{\alpha} = 2 \mathcal{P}_{nn+1} - \mathbb{I}_{nn+1}.
$$
An analogous situation arises in the anisotropic $XXZ$ chain
\begin{equation}
\label{HXXZ}
H_{XXZ} = \sum _{n=1}^{N-1} \left( \sigma_n^x \sigma_{n+1}^x + \sigma_n^y \sigma_{n+1}^y  + \Delta\sigma_n^z \sigma_{n+1}^z  \right) + \frac{\left( q - q ^{-1} \right)}{2} \left( \sigma_1^z - \sigma_{N}^z \right) ,
\end{equation}
which commutes \c{PasSal} with the global generators of the quantum algebra $\mathcal{U}_q(\mathfrak{sl}(2))$ \c{KR}. Here the role of the second algebra is played by the Temperley-Lieb algebra $TL_N(q)$, whose
generators in the space $\mathcal{H}= \underset {1}{\overset {N}{\otimes}} \mathbb{C}^2$ coincide with the constant $R$-matrix ($\omega(q)=q-1/q$)
\begin{equation}
\label{RcXXZ}
\check{R} _{XXZ} (q) = 
\left(\begin{array}{cccc}
q & 0 & 0 & 0 \\
0 & \omega(q) & 1 & 0 \\
0 & 1& 0 & 0 \\
0 & 0 & 0 & q 
\end{array} \right).
\end{equation}
As in the case of the $XXX$ spin chain, the Hamiltonian \Ref{HXXZ} can be expressed in terms of the
generator \Ref{RcXXZ} of the algebra $TL_N(q)$ ($\Delta=(q+q^{-1})/2$)
\begin{eqnarray*}
&&\sigma_n^x \sigma_{n+1}^x + \sigma_n^y \sigma_{n+1}^y  + \Delta\sigma_n^z \sigma_{n+1}^z  + \frac{\omega(q)}{2} \left( \sigma_n^z - \sigma_{n+1}^z \right) \nonumber \\
&& = 2 \check{R} _{XXZ} (q) - \left( \frac{\omega(q)}{2} + q \right) \mathbb{I}_{nn+1} . \nonumber
\end{eqnarray*}
The dynamical symmetry algebra of the $XXZ$ chain is the quantum 
affine algebra $\mathcal{U}_q(\widehat{\mathfrak{sl}_2})$ \c{JimboMiwa}.

The Temperley-Lieb algebra is a quotient of the Hecke algebra (see section 2) and allows for an $R$-matrix representation in the space $\mathcal{H}= \underset {1}{\overset {N}{\otimes}} \mathbb{C}^n$ for any $n=2,3,\ldots$. There is corresponding spectral parameter depending R-matrix obtained by the
Yang - Baxterization process. Consequently, it is possible to construct an integrable spin chain \c{Kulish}. The open spin chain Hamiltonian is the sum of the $TL_N(q)$ generators $X_j$
$$
H_{TL}= \sum X_j,
$$
where $X_j$ act nontrivially on $\mathbb{C} _j^n \otimes \mathbb{C}_{j+1}^n$ and as the identity  matrix on the other factors of $\mathcal{H}$. The aim of this work is to describe the quantum algebra $\mathcal{U}_q(n)$ which is the symmetry algebra of such spin system and to show that the structures (categories) of finite dimensional representations of these algebras $\mathcal{U}_q(n)$ and $\mathfrak{sl}_2$ coincide. In this case $\mathcal{U}_q(n)$ and $TL_N(q)$ are each other's centralizers in the space $\mathcal{H}= \underset {1}{\overset {N}{\otimes}} \mathbb{C}^n$. We consider the general case when the complex parameter $q \in \mathbb{C}^{\ast}$ is not a root of unity.

Let us note that the relation between $TL_N(q)$ and integrable spin chains was actively used in many works and monographs (see for example \c{Baxter,BRJ,DWA,Martin,BK, KoberleLS,MR,GNPR,Doikou} and the references within). However, the authors used particular realizations of the generators $X_j$, related to some Lie algebras (or quantum algerbas). Characteristic property of the latter ones was
the existence of one-dimensional representation in the decomposition of the tensor product of two fundamental representations $V_j \otimes V_{j+1}$. Then $X_j$ was proportional to the rank one projector on this subspace, and the symmetry algerba was identified with the choosen algebra. We point out that the symmetry algerba $\mathcal{U}_q(n)$ is bigger and its Clebsch - Gordan decomposition of $V_j \otimes V_{j+1}$ has only two summands similar to the sl(2) case $\mathbb{C}^2 \otimes \mathbb{C}^2 = \mathbb{C}^3 \oplus \mathbb{C}^1$.

The dual Hopf algebra $\mathcal{U}_q(n)^{\ast}$ was introduced 
as the quantum group of nondegenerate bilinear form in \cite{DVL,Gu}. The categories of co-modules of $\mathcal{U}_q(n)^{\ast}$ and their generalisations were studied in \c{Bichon, EO} where it was shown that the 
categories of co-modules of $\mathcal{U}_q(n)^{\ast}$ are equivalent 
to the category of co-modules of the quantum group $SL_q(2)$.

This paper is organised as follows. In Section 2 we recall the definition of  Hecke algebra and of its useful quotient,  the Temperley-Lieb algebra. We also define a representation of the $TL_N(q)$ algebra using an arbitrary invertible $n\times n$ matrix. In Section 3 we use this representation and the $R$-matrix derived from it to define the quantum algebra $\mathcal{U}_q(n)$. We show that the algebras $\mathcal{U}_q(n)$ and $TL_N(q)$ are each other's centralizers in the space $\mathcal{H}=\underset {1}{\overset {N}{\otimes}} \mathbb{C}^n$. We discuss the decomposition of tensor products of irreducible representations of $\mathcal{U}_q(n)$. It turns out that a representation ring of the algebra
$\mathcal{U}_q$ is equivalent to the representation ring of the $\mathfrak{sl}_2$ Lie algebra. 
%
%

\section{Hecke and Temperley-Lieb Algebras}

Both  algebras $\mathcal{H}_N(q)$ and $TL_N(q)$ are quotients of the group algebra of the braid group  $\mathcal{B}_N$ generated by $(N-1)$ generators $\check{R}_j$, $j=1,2,\ldots,N-1$, their inverses $\check{R}_j^{-1}$ and the relations (see \c{CP})
\begin{eqnarray}
\check{R}_j\check{R}_k & = & \check{R}_k \check{R}_j , \ |j-k| > 1, \nonumber \\
\label{BG}
\check{R}_j\check{R}_k  \check{R}_j & = & \check{R}_k \check{R}_j  \check{R}_k ,   \ |j-k| =1.
\end{eqnarray}
The Hecke algebra $\mathcal{H}_N(q)$ is obtained by adding to these 
relations the following characteristic equations obeyed by generators
\begin{equation}
\label{cheqRj}
\left( \check{R}_j - q \right) \left( \check{R}_j + 1/q \right) = 0. 
\end{equation}
It is known that $\mathcal{H}_N(q)$ is isomorphic to the group algebra $\mathbb{C}[S_N]$. 
Consequently, irreducible representations of the Hecke algebra, as of $S_N$, are parametrized by Young diagrams. By virtue of \Ref{cheqRj} we can write $\check{R}$ using the idempotents $P_+$ and $P_-$ ($P_++P_-=1$):
\begin{equation}
\label{Rproj}
\check{R} = q P_+ - \frac{1}{q}P_- = q \mathbb{I} - \left( q + \frac{1}{q} \right) P_- := q \mathbb{I} + X
\end{equation}
Substituting the expression \Ref{Rproj} for $\check{R}$ in terms of $X$, into the braid group relations \Ref{BG} one gets relations for $X_j, X_k$, $|j-k|=1$
\begin{equation}
\label{BGX}
X_jX_kX_j - X_j = X_kX_jX_k - X_k . 
\end{equation}
Requiring that each side of \Ref{BGX} is zero we obtain the quotient algebra of the Hecke algebra, the Temperley-Lieb algebra $TL_N(q)$.  It is defined by the generators $X_j$, $j=1,2,\ldots, N-1$ and the relations ($\nu(q)=q+1/q$):
\begin{eqnarray}
&&X_j^2 =  - \left( q + \frac{1}{q} \right) X_j = - \nu (q) X_j, \nonumber \\
\label{TLX}
&&X_jX_kX_j = X_j , \quad | j - k | = 1. 
\end{eqnarray}
The dimension of the Hecke algebra is $N!$, the same as the dimension 
of the symmetric group $S_N$,  the dimension of $TL_N(q)$ is 
equal to the Catalan number $C_N=(2N)!/N!(N+1)!$. In connection with integrable spin systems we will be interested in representations of $TL_N(q)$ on the tensor product space $\mathcal{H}=\underset {1}{\overset {N}{\otimes}} \mathbb{C}^n$. One representation is defined by an invertible $n \times n$ matrix 
$b \in GL(n,\mathbb{C})$ which can also be seen as an $n^2$ dimensional vector 
$\{ b_{cd} \} \in \mathbb{C}^n \otimes \mathbb{C}^n$ \c{MR}. We use the notation 
$\bar{b}:=b^{-1}$ and view this matrix also as an $n^2$ dimensional vector 
$\{\bar{b}_{cd} \} \in \mathbb{C}^n \otimes \mathbb{C}^n$. The generators $X_j$ can be expressed as 
\begin{equation}
\label{Xbb}
\left( X_j \right)_{cd,xy}  = b_{cd} \bar{b}_{xy} \in \mathrm{Mat} \left( \mathbb{C}_j^n \otimes \mathbb{C}_{j+1}^n \right),
\end{equation}
where we explicitely write the indices corresponding to the factors in the tensor product space $\mathcal{H}$. It is easy to see, that the second relation \Ref{TLX} is automatically satisfied and the first one determines the parameter $q$ ($\nu(q)=q+1/q$):
\begin{equation}
\label{Xjbb}
X_j^2 = X_j  \, \mathrm{tr} \, b^t\bar{b}, \quad \mathrm{tr} \, b^t \bar{b}= - \left( q + \frac{1}{q} \right) = - \nu (q) .
\end{equation}
An obvious invariance of the braid group relations (the Yang - Baxter equation) \Ref{BG} in this
representation with respect to the transformation of the R-matrix
$$
 \check{R} \to  \mathrm{AdM} \otimes \mathrm{Ad M} \left(\check{R}\right) ,  \qquad \mathrm{M} \in \mathrm{GL}(n,\mathbb{C}) ,
$$
results in the following transformation of the matrix 
$$
b \to \mathrm{M} b \mathrm{M}^t.
$$ 
If one uses an $R$-matrix depending on a spectral parameter 
(Yang - Baxterization of $\check{R}(q)$)
\begin{equation}
\label{RBax}
\check{R}(u;q)= u \check{R}(q) - \frac{1}{u}\left( \check{R}(q) \right)^{-1} = \omega(uq) \mathbb{I} + \omega (u)X,
\end{equation}
where $ \check{R}(q) ^{-1} = (1/q)\mathbb{I} + X$, then relation \Ref{TLX} can be written as
\begin{equation}
\label{TLR}
\check{R}(q^{-1};q)\check{R}(q^{-2};q)\check{R}(q^{-1};q)=0.
\end{equation}
In terms of constant $R$-matrices (generators of $TL_N(q)$) this relation has the form\break ($|i-k|=1$, $\check{R}_i=\check{R}_{i,i+1}$)
\begin{equation}
\label{TLRc}
\left( \check{R}_i - q \mathbb{I} \right)
\left( \nu(q) \check{R}_k - q^2 \mathbb{I} \right) 
\left( \check{R}_i - q \mathbb{I} \right) = 0
\end{equation}
Replacing in \Ref{TLR} the expression $\check{R}(u;q)= \omega(u) \check{R}(q)+u^{-1}\omega(q)\mathbb{I},$ or in \Ref{TLX} substituting $X=\check{R}-q I$ yields the vanishing of the $q$-antisymmetriser
\begin{equation}
\label{ }
\mathbb{I} - q^{-1} \left( \check{R}_{12} + \check{R}_{23} \right) + q^{-2} \left( \check{R}_{12} \check{R}_{23} + \check{R}_{23}\check{R}_{12} \right) - q^{-1} \check{R}_{12} \check{R}_{23} \check{R}_{12} = 0.
\end{equation}

Thus the irreducible representations of $TL_N(q)$ are parametrized by Young diagrams containing only two rows.

The constructed representation \Ref{Rproj}, \Ref{Xbb} is 
reducible. The decomposition of this representation into the 
irreducible ones  will be discussed in the next section.

%
%

\section{Quantum Algebra $\mathcal{U}_q(n)$}

According to the $R$-matrix approach to the theory of quantum groups \c{RTF}, the R-matrix defines relations between the generators of the quantum algebra $\mathcal{U}_q$ and its dual Hopf algebra, the quantum group $\mathcal{A} (R)$. In this paper the emphasis will be on the 
quantum algebra $\mathcal{U}_q$ and its finite dimensional representations $V_k$, $k=0, 1, 2, \dots $. 
The generators of  $\mathcal{U}_q$ can be identified with the  $L$-operator ($L$-matrix) entries and their exchange relations (commutation relations) 
follow from the analogue of the Yang-Baxter relations
\begin{equation}
\label{YBR}
\check{R} _{12} L _{a_1q} L _{a_2q}= L _{a_1q}L _{a_2q}\check{R} _{12} ,
\end{equation}
where the indices $a_1$ and $a_2$ refer to the representation 
spaces $V_{a_1}$ and $V_{a_2}$, respectively, and index $q$ 
refers to the algebra $\mathcal{U}_q$. Hence the equation \Ref{YBR} is given in $\mathrm{End}\left(V_{a_1} \otimes V_{a_2} \right) \otimes \mathcal{U}_q$.

In general the $L$-operator is defined through the universal $R$-matrix, where a finite dimensional representation is applied to one of the factors of the universal $R$-matrix 
\begin{equation}
\label{RUn}
\mathcal{R}_{univ} = \sum _j \mathcal{R}_1^{(j)} \otimes \mathcal{R}_2^{(j)} := 
\mathcal{R}_1 \otimes \mathcal{R}_2 \in \mathcal{U}_q \otimes \mathcal{U}_q , 
\end{equation}
with 
\begin{equation}
\label{LMat}
 L _{aq} = \left( \rho \otimes \mathrm{id} \right) 
 \mathcal{R}_{univ} = 
 \rho(\mathcal{R}_1)\otimes \mathcal{R}_2, 
\end{equation}
where $\rho : \mathcal{U}_q \to \mathrm{End}(V_a)$. 

Furthermore, the universal $R$-matrix satisfies Drinfeld's axioms of the quasi-triangular Hopf algebras
 \c{Drinfeld, CP}. 
In particular,
\begin{equation}
\label{ }
\left( \mathrm{id} \otimes \Delta \right) \mathcal{R}= \mathcal{R}_{13} \mathcal{R}_{12} . 
\end{equation}
Thus, choosing the appropriate representation space as the first space, one obtains the co-product of the generators of $\mathcal{U}_q$ from the following matrix equation
\begin{equation}
\label{Copr}
\left( \mathrm{id} \otimes \Delta \right) L _{aq} = L _{aq_2} L _{aq_1} \in \mathrm{End}(V_a) \otimes  \mathcal{U}_q \otimes \mathcal{U}_q .
\end{equation}
The case when $V_a = \mathbb{C}^3$ is of particular interest 
and it will be studied below in detail. To this end the generators of $\mathcal{U}_q$ are denoted by $\{ A_i , B_i , C_i , i = 1, 2, 3 \}$ and the $L$-matrix
is given by
\begin{equation}
\label{L-mat}
L _{aq} = \left(\begin{array}{ccc}
A_1 & B_1 & B_3 \\
C_1 & A_2 & B_2 \\
C_3 & C_2 & A_3 \end{array} \right).
\end{equation}
Multiplying two $L$-matrices with entries in the corresponding factors $\mathcal{U}_q(n) \otimes \mathcal{U}_q(n)$ we obtain
\begin{equation}
\label{LCopr}
\Delta \left( L_{ab} \right) = \sum _{k=1}^3  L_{kb}  \otimes  
L_{ak}= \sum _{k=1}^3 \left( \mathbb{I} \otimes L_{ak}\right)
\left( L_{kb} \otimes \mathbb{I} \right) ,
\end{equation}
or explicitly for the generators
\begin{eqnarray}
\label{B1Copr}
\Delta (B_1) & = & B_1 \otimes A_1 + A_2 \otimes B_1+ C_2 \otimes B_3,  \\
\label{B2Copr}
\Delta (B_2) & = & B_3 \otimes C_1 + B_2 \otimes A_2+ A_3 \otimes B_2,  \\
\label{B3Copr}
\Delta (B_3) & = & B_3 \otimes A_1 + B_2 \otimes B_1+ A_3 \otimes B_3 ,
\end{eqnarray}
etc. The central element in $\mathcal{U}_q$ is obtained from the defining relation \Ref{YBR}
\begin{equation}
\label{Casimir}
b^{-1} L _{aq} b L _{aq}^t = c_2 \mathbb{I},
\end{equation} 
$$
c_2 = \sum_{jkl} \left( b^{-1} \right)_{1j}L_{jk}b_{kl}L_{1l} .
$$
%
However this central element is group-like: $\Delta c_2 = c_2 \otimes c_2$. 
It is proportional to the identity in the tensor product of representations. 
The analogue of the $\mathcal{U}_q(2)$ Casimir operator can be 
obtained according to \c{RTF} using $L_+:= L$ and $L_-:=
\left( \rho \otimes \mathrm{id} \right) 
\left(\mathcal{R}_{21} \right)^{-1}$ as 
$\mathrm{tr}_q L_+L_-^{-1} = \mathrm{tr} \, b\bar{b}^tL_+L_-^{-1}$.  

In the case when $V_a$, $V_q $ are the three dimensional space $V_a\simeq V_q \simeq \mathbb{C}^3$ and the $b$ matrix is taken from the references \c{BK,KoberleLS} 
\begin{equation}
\label{bmat}
b_{ij} = p ^{2-i}\delta _{i4-j}=\left( b^{-1}\right)_{ij}
\end{equation}
$c_2$ is written as 
$$
c_2=p\left( \frac{1}{p} A_3A_1 + C_2B_1 + p C_3B_3\right).
$$
Parameters $p$ and $q$ are related: $p^2+1+p^{-2}=-(q+q^{-1})$. For the explicit $L$-operator and its
$3\times 3$ blocks we get $c_2=q I_q$ where $I_q$ is the identity operator on $V_q=\mathbb{C}^3$. The form of the generators $\{A_i,B_i,C_i, i=1,2,3\}$ which corresponds to the choice of the $b$ matrix \Ref{bmat}, 
follows from the expression for the $\check{R}$ and $L$-matrices
\begin{equation}
\label{Xgen}
X = |b\rangle \otimes \langle b| = 
\left(\begin{array}{ccc|ccc|ccc}
0 &  &  &  &  &  &  &  &  \\ 
& 0 &  &  &  &  &  &  &  \\ 
&  & p^2 &  & p &  & 1 &  &  \\\hline 
&  &  & 0 &  &  &  &  &  \\ 
&  & p &  & 1 &  & 1/p &  &  \\ 
&  &  &  &  & 0 &  &  &  \\\hline 
&  & 1 &  & 1/p &  & 1/p^2 &  &  \\ 
&  &  &  &  &  &  & 0 &  \\ 
&  &  &  &  &  &  &  & 0\end{array}\right),
\end{equation}
and
\begin{eqnarray}
& & L = \mathcal{P} \check{R} =   \mathcal{P} \left( qI + X \right) \nonumber \\
\label{Rmat}
& & = \left(\begin{array}{ccc|ccc|ccc}
q &  &  &  &  &  &  &  &  \\ 
   & 0 &  & q &  &  &  &  &  \\ 
   &  & 1 &  & p^{-1} &  & q+p^{-2} &  &  \\\hline
   & q &  & 0 &  &  &  &  &  \\ 
   &  & p &  & q+1 &  & p^{-1} &  &  \\ 
   &  &  &  &  & 0 &  & q &  \\ \hline
   &  & q+p^2 &  & p &  & 1 &  &  \\ 
   &  &  &  &  & q &  & 0 &  \\ 
   &  &  &  &  &  &  &  & q\end{array}\right) .
\end{eqnarray}
where $\mathcal{P}$ is the permutation matrix. For example, we have
\begin{equation}
\label{Bops}
B_1 = \left(\begin{array}{ccc}0 & 0 & 0 \\q & 0 & 0 \\0 & p^{-1} & 0\end{array}\right), \
B_2 = \left(\begin{array}{ccc}0 & 0 & 0 \\p^{-1} & 0 & 0 \\0 & q & 0\end{array}\right).
\end{equation}

If we choose $b$ as in equation \Ref{bmat} the $R$-matrix \Ref{Rmat} commutes with $h\otimes 1+1 \otimes h$ where $h=\mathrm{diag}(1,0,-1)$. 
As a highest spin vector (pseudovacuum of the corresponding integrable spin chain \c{Fadddeev}) we choose $\theta=(1,0,0)^t \in \mathbb{C}^3$ \c{Kulish, KoberleLS}. 
If we act on the tensor product of these vectors $\theta \otimes \theta \in \mathbb{C}^3\otimes \mathbb{C}^3$ with the coproduct of the lowering operators $\Delta(B_j)$, $j=1,2,3$ we obtain new vectors. By looking at the explicit forms of the operators $A_j$,$B_j$,$C_j$ in the space $\mathbb{C}^3$ we can convince ourselves that the vectors $(\Delta(B_i))^k \theta \otimes \theta$, $i=1,2; k=1,2,3$ are linearly independent. 
Together with the vectors $\theta\otimes \theta$ and
$$
(\Delta(B_i))^4 \theta \otimes \theta \simeq \left(\begin{array}{c}0 \\0 \\1\end{array}\right) \otimes \left(\begin{array}{c}0 \\0 \\1\end{array}\right),
$$
they span an $8$ dimensional subspace. The vector $\Delta(B_3)\theta \otimes \theta$ is a linear combination of the vectors $(\Delta(B_i))^2 \theta \otimes \theta$, $i=1,2$. The vector $|b\rangle=(0 0 p|010| p^{-1} 0 0)^t$ spans a one dimensional invariant subspace. Thus we have the following decomposition
\begin{equation}
\label{n=2}
\mathbb{C}^3 \otimes \mathbb{C}^3 = \mathbb{C}^8 \oplus \mathbb{C} . 
\end{equation}
This decomposition can also be obtained using the projectors $P_+$, $P_-$ \Ref{Rproj}, \Ref{Xbb} 
expressed in terms of $b$ matrix (vector) \Ref{bmat}. 
Due to the commutativity of the $R$-matrix 
$\check{R}_{q_1q_2}$ with the co-product \Ref{Copr}, \Ref{LCopr} the corresponding 
subspaces $P_{\pm} (\mathbb{C}^3 \otimes \mathbb{C}^3)$ are invariant. 

Similarly, using $\Delta^3(Y), Y \in \mathcal{U}_q(3)$ one can get the 
decomposition of 
\begin{equation}
\label{n=3}
\mathbb{C}^3 \otimes \mathbb{C}^3 \otimes \mathbb{C}^3 = 
\mathbb{C}^{21} \oplus \mathbb{C}^3 \oplus \mathbb{C}^3. 
\end{equation}
This type of decomposition is valid for any n.

The result of applying the co-product $\Delta$, given by \Ref{Copr}, on the generators of the quantum algebra $\mathcal{U}_q$ several times can also be presented in the matrix form
\begin{equation}
\label{TN}
\left( \mathrm{id} \otimes \Delta^N \right) L _{aq} = L _{aq_N}  \ldots L _{aq_2} L _{aq_1} := T ^{(N)},
\end{equation}
where $\Delta^N : \mathcal{U}_q \to \left( \mathcal{U}_q\right)^{\otimes N}$, $\Delta^1:= \mathrm{id}$,
$\Delta^2:= \Delta$, $\Delta^3:=\left( \mathrm{id} \otimes \Delta \right) \circ \Delta$, etc. In general case of the tensor representation of $TL_N(q)$, with the space $\mathbb{C}^n$ at each site, the generators of the algebra $\check{R}_{q_k q_{k+1}}:= \check{R}_{k,k+1}$ commute with the generators \Ref{TN} of the global (diagonal) action of the quantum algebra $\mathcal{U}_q(n)$ in the space $\mathcal{H} =  \underset {1}{\overset {N}{\otimes}} \mathbb{C}^n$. This follows from the relation
$$  
\check{R}_{k,k+1} L_{aq_{k+1}} L_{aq_k} = L_{aq_{k+1}} L_{aq_k} \check R_{k,k+1},
$$
and the possibility due to the co-associativity of the coproduct
to write the product of of $L_{aq_j}$ \Ref{Rmat} as
$$ T^N = L_{aq_N} \ldots L_{aq_{k+2}} \Delta_k(L_{aq_k}) L_{aq_{k-1}}  \ldots L_{aq_1}.$$
Hence, 
$$
\check{R}_{k,k+1} T ^{(N)} = T ^{(N)}  \check{R}_{k,k+1} . 
$$
Thus, the algebras $\mathcal{U}_q(n)$ and $TL_N(q)$ are each other's centralizers in the space $\mathcal{H}$. The tensor representation of $TL_N(q)$ in 
$\mathcal{H}$ decomposes into irreducible factors whose multiplicities 
are given by the dimensions of the irreducible representations of 
the algebra $\mathcal{U}_q(n)$, corresponding to the same Young diagrams
\begin{equation}
\label{decom}
\mathcal{H} = \underset {k}{\overset {N}{\otimes}} \mathbb{C}^n = \underset {k}{\overset {N}{\oplus}} 
p_k(n) W_k(N) =  \underset {k}{\overset {N}{\oplus}} \nu_k(N) V_k(n).
\end{equation}
In this decomposition the index $k$ parametrizes the Young diagrams with two rows and $N$ boxes and multiplicities are given by the dimensions of the corresponding irreducible representations
\begin{equation}
\label{ }
p_k(n) = \dim V_k(n), \  \nu_k(N) = \dim W_k(N) ,
\end{equation}   
of the algebras $\mathcal{U}_q(n)$ and $TL_N(q)$, respectively. 
As for the finite dimensional irreducible  representations of the Lie algebra $\mathfrak{sl}_2$, $V_0(n)= \mathbb{C}$ is the one-dimensional (scalar) representation and the fundamental representation of the algebra $\mathcal{U}_q(n)$ is  $n$ dimensional, $V_1(n) \simeq\mathbb{C}^n$. The dimensions of other representations follow from the trivial multiplicities of the factors in the decomposition of the tensor product of the $V_k(n)$ and the fundamental representation $V_1(n)$ into two irreducible factors, as for the $\mathfrak{sl}_2$, 
\begin{equation}
\label{ }
V_1(n) \otimes V_k(n) = V_{k+1}(n) \oplus V_{k-1}(n) .
\end{equation}
Thus, for the dimensions $p_k(n) = \dim V_k(n)$ the following recurrence relation is valid 
\begin{equation}
\label{ }
n \cdot p_k(n) = p_{k+1}(n) + p_{k-1}(n), 
\end{equation}
with the initial conditions $p_{-1}(n)=0$, $p_0(n)=1$, whose solutions are Chebyshev polynomials of the second kind
\begin{equation}
\label{ }
p_k(n) = \frac{\sin (k+1) \theta}{\sin \theta} , \ n = 2 \cos \theta.
\end{equation} 

The multiplicities $\nu_k(N)$, or the dimensions of the subspaces $W_k(N)$ in \Ref{decom} 
is the number of paths that go from the top of the Bratteli diagram to the Young diagram corresponding to the representation $W_k(N)$. If $\lambda \vdash N$ is the partition of $N$, $\lambda = (\lambda_1 \geqslant \lambda _2 | \lambda_1 + \lambda _2 = N )$, then $k = \lambda_1 - \lambda _2$ and
\begin{equation}
\label{ }
\nu_k(N) = \nu_{k+1} (N-1) + \nu_{k-1}(N-1).
\end{equation}

The subspaces invariant under the diagonal action of the quantum algebra 
$\mathcal{U}_q(n)$ on the space $\mathcal{H}$, can be obtained using the projectors (idempotents), which can be expressed in terms of the elements of the Temperley-Lieb algebra $TL_N(q)$. Using the $R$-matrix depending on a spectral parameter, the projector $P_N^{(+)}$ on the symmetric subspace can be written in the following way \c{KS,Isaev}
\begin{equation}
\label{ }
P_N^{(+)} \simeq P_{N-1}^{(+)} \check{R}_{N-1N}(q^{N-1}; q)  P_{N-1}^{(+)},
\end{equation} 
where $P_2^{(+)} = \mathbb{I} - P_-, P_-=- X/\nu(q)$. 

It is possible also to construct quantum spaces which are covariant with respect to the action (respectively co-action) of the algebra $\mathcal{U}_q(n)$ (respectively the dual algebra $\mathcal{U}_q(n)^{\ast}$ \c{Bichon}). To this end, let us consider the tensor algebra $T(V)$ over the fundamental representation space $V \simeq \mathbb{C}^n$ and its quotients by the corresponding ideals. If one denotes the basis elements (generators) of $V$ by $Z=(z_1,z_2,\ldots,z_n)^t$ and chooses the relations which generate the ideal $J$ ($Z_1Z_2=Z\otimes Z$) as:
\begin{equation}
\label{ }
\check{R}Z_1Z_2=q Z_1Z_2, \  \mathrm{or} \ X Z_1Z_2=0,
\end{equation}
then one gets the analogue of the symmetric subalgebra $S(V)\simeq T(V)/J$. The ideal $J$ is also generated by a single quadratic relation $Z^t\bar{b}Z=0$. For the analogue of the antisymmetric subalgebra $\Lambda(V)\simeq T(V)/J_-$ we choose relations (see \Ref{Rproj})
\begin{equation}
\label{ }
\check{R} W_1 W_2 = -\frac{1}{q}W_1 W_2, \ \mathrm{or} \ (I-P_-)W_1 W_2=0,
\end{equation} 
where we denoted the generators by $W=(w_1,w_2,\ldots,w_n)^t$. Among quadratic combinations there is only one different from zero: $W^t \bar{b} W=\sum w_{\alpha} (b^{-1})_{\alpha \beta} w_{\beta}$ and the Hilbert-Poincar\'e series of $\Lambda(V)$ is equal to $P_-(t)=1+nt+t^2$. One can recognize the generating function of the Chebyshev polynomials, and  the relation $P_+(t) = 1/P_-(-t) = \sum p_k(n) t^k$ reproduces the dimensions of the homogeneous subspaces of $S(V)$.

%
%

\section{Acknowledgements}

We acknowledge useful discussions with P. ~Etingof and A. ~Mudrov. 
This work was supported by the grants RFFI 05-01-00922, NS-5403.2006.1 
and the FCT project POCTI/MAT/58452/2004. In addition to that Z. Nagy 
benefited from the FCT grant SFRH/BPD/25310/2005. Part of this work was done during a visit to the Laboratoire de Physique Th\'eorique et Mod\'elisation of the University of Cergy-Pontoise.

\textbf{Note added in proofs:} We are grateful to the referee for calling our attention to the references \c{DM} and \c{RS}.


\end{document}